\documentclass[leqno,12pt]{amsart}
\usepackage{latexsym}
\usepackage{mathrsfs}
\usepackage{amsmath}
\usepackage{amssymb}
\usepackage[bookmarks]{hyperref}
\usepackage{pstricks,pst-plot}

\newtheorem{theorem}{Theorem}[section]
\newtheorem{lemma}[theorem]{Lemma}
\newtheorem{corollary}[theorem]{Corollary}
\newtheorem{proposition}[theorem]{Proposition}

\theoremstyle{definition}

\newtheorem{definition}[theorem]{Definition}

\newcommand{\B}{\mathbb{B}}
\newcommand{\C}{\mathbb{C}}

\newcommand{\Q}{\mathbb{Q}}

\newcommand{\Z}{\mathbb{Z}}
\renewcommand{\P}{\mathbb{P}}
\newcommand{\R}{\mathbb{R}}

\newcommand{\bthm}{\begin{theorem}}
\newcommand{\ethm}{\end{theorem}}
\newcommand{\blem}{\begin{lemma}}
\newcommand{\elem}{\end{lemma}}
\newcommand{\bcor}{\begin{corollary}}
\newcommand{\ecor}{\end{corollary}}
\newcommand{\bprop}{\begin{proposition}}
\newcommand{\eprop}{\end{proposition}}
\newcommand{\bdefn}{\begin{definition}}
\newcommand{\edefn}{\end{definition}}
\newcommand{\bpf}{\begin{proof}}
\newcommand{\epf}{\end{proof}}

\def\hr#1{h_r({#1})}
\def\hhr#1{h_r({#1}) \sm {#1}}

\def\vep {\varepsilon}
\def \sm {\setminus}
\def \ob {\overline \B}

\def\h#1{\widehat {#1}}
\def\hh#1{\widehat {#1} \sm  {#1}}

\def\ol {\overline}

\def \ma {\mathfrak{M}_A}  
\def \mb {\mathfrak{M}_B}

\def\pb{\partial B}

\def\P{{\mathscr P}}
\def\Q{{\mathscr Q}}
\def\vector#1#2{(#1_1,#1_{#2})}
\def\row#1#2{#1_1,\ldots,#1_{#2}}

\begin{document} \bibliographystyle{plain}
\title[Polynomially convex sets whose union has hull]{Polynomially convex sets\\ whose union has nontrivial hull}
\author{Alexander J. Izzo}
\address{Department of Mathematics and Statistics, Bowling Green State University, Bowling Green, OH 43403}
\email{aizzo@bgsu.edu}
\thanks{The author was supported by NSF Grant DMS-1856010 and by a Simons Foundation Collaboration Grant.}


\begin{abstract}
Several results 
$\vphantom{\widehat{\widehat{\widehat{\widehat{\widehat{\widehat{\widehat X}}}}}}}$concerning pairs of polynomially convex sets whose union is not even rationally convex are given.  It is shown that there is no
restriction on how two spaces can be embedded in some $\C^N$ so as to be polynomially convex but have nonrationally convex union.  It is shown that there exist two disjoint polynomially convex Cantor sets in $\C^3$ whose union is not rationally convex.  The analogous assertion for arcs is also established.  As an application it is shown that every simple closed curve in $\C^N$, $N\geq 3$, can be approximated uniformly by locally polynomially convex simple closed curves that are not rationally convex.
\end{abstract}

\maketitle 

\vskip -2.40 true in
\centerline{\footnotesize\it Dedicated to the memory of Peter Duren} 
\vskip 2.40 truein

\section{Introduction}

The motivation for the work presented here is a recent result of the author and Lee Stout \cite[Theorem~3.1]{Izzo-Stout:2022}: Given an arbitrary polynomially convex simple closed curve $\gamma$, every rectifiable simple closed curve $\sigma$
sufficiently close to 
$\gamma$ in the uniform metric is also polynomially convex.  Note that although there is no regularity hypothesis on the curve $\gamma$, the curve $\sigma$ is required to be rectifiable.  In fact, it follows from the existence of nonpolynomially convex arcs that the statement becomes false without the rectifiability hypothesis; one can obtain a nonpolynomially convex simple closed curve arbitrarily close to $\gamma$ by choosing a nonpolynomially convex arc contained in a small neighborhood of some point of $\gamma$ and modifying $\gamma$ so as to contain that nonpolynomially convex arc.  (The existence of a nonpolynomially convex arc was first proved by John Wermer \cite{Wermer:1955} and can be found in standard texts, e.g., \cite[pp.~222-223]{Alexander-Wermer:1998}, \cite[pp.~29 and 70]{Gamelin:1984}, and \cite[pp.~53--54]{Stout:2007}.)   However, every \emph{rectifiable} arc is polynomially convex, and hence every \emph{rectifiable} simple closed curve is \emph{locally} polynomially convex.  One might, therefore, conjecture that the above result of the author and Stout could be generalized to assert that every locally polynomially convex simple closed curve sufficiently near a polynomially convex simple closed curve is itself polynomially convex.  One of the main goals of the present paper is to show that, on the contrary, every simple closed curve in $\C^N$, $N\geq3$, can be approximated uniformly by a locally polynomially convex, nonpolynomially convex simple closed curve.  In fact, the approximating simple closed curve $\sigma$ can be taken to be the union of \emph{two} polynomially convex arcs whose interiors cover $\sigma$ (and such that, in fact, $\sigma$ is not even rationally convex) (Theorem~\ref{path-curve}).  Whether the same holds in $\C^2$ remains open.

One possible approach to trying to prove the assertion just made regarding curves in $\C^N$, $N\geq 3$, is to attempt to construct two disjoint arcs near $\gamma$ each of which is polynomially convex but whose union is not polynomially (or rationally) convex, and then to attach to these two disjoint arcs additional arcs to form a simple closed curve.  We are thus led to the question of whether there exist two disjoint polynomially convex arcs whose union is not polynomially convex.  More generally, the question arises as to what disjoint polynomially convex sets can have nonpolynomially convex union.  A fundamental question here is whether there exist two disjoint polynomially convex Cantor sets whose union is not polynomially convex.

We will establish quite general results regarding the issue of pairs of polynomially convex sets whose union is not polynomially convex.  Roughly, we will show that there is \emph{no} restriction on how two spaces can be embedded in some $\C^N$ so as to be polynomially convex but have nonpolynomially convex union, and even nonrationally convex union (Theorem~\ref{general-unions} and Corollary~\ref{n-dimensional-case}) .  We will also obtain several results about embedding particular spaces in $\C^N$ for particular values of $N$.  In particular, we will show that there exist two disjoint polynomially convex Cantor sets in $\C^3$ whose union is not rationally convex (Corollary~\ref{disjoint-Cantor-sets}).  The same holds for arcs in place of Cantor sets (Corollary~\ref{arc-union}).  We will also obtain two polynomially convex Cantor sets in $\C^2$ whose union is not rationally convex (Theorem~\ref{two-Cantor-sets}); however, in that setting we do not know whether the polynomially convex Cantor sets can be taken to be disjoint.

An interesting related question is whether there exists a minimal nonpolynomially convex Cantor set in some $\C^N$, that is, whether there is a nonpolynomially convex (or nonrationally convex) Cantor set every closed subset of which is polynomially convex.  We will not answer this question for sets in a finite-dimensional $\C^N$, but we will answer the question affirmatively in a more general abstract uniform algebra setting.

Stout raised the question of whether the union of the closed unit ball in $\C^N$ and a polynomially convex arc that meets the ball precisely in a single end point must be polynomially convex.  When the arc is rectifiable, polynomial convexity of the union is a consequence of a standard result \cite[Theorem~3.1.1]{Stout:2007}.  As noted by Stout, applying \cite[Theorem~1.6]{Izzo-Stout:2022b} shows that polynomial convexity of the union holds more generally whenever the interior of the arc is locally rectifable (Theorem~\ref{Stout-observation}).  We will show that nevertheless, in general the union of the closed unit ball and a polynomially convex arc that meets the ball in precisely a single end point need not be even rationally convex (Corollary~\ref{non-rat-ball-n-chain}).  We will also give an example of a polynomially convex arc disjoint from the closed unit ball such that the union is not rationally convex.

It was once conjectured that whenever the polynomial hull $\h X$ of a compact set $X$ in $\C^N$ is strictly larger than $X$, the complementary set $\hh X$ must contain an analytic disc.  This conjecture was disproved by Gabriel Stolzenberg \cite{Stolzenberg:1963}.  There are now many known counterexamples due to several mathematicians.  By combining our results on pairs of polynomially convex sets whose union has nontrivial hull with a construction of nontrivial hulls without analytic discs in the paper of the author and Norman Levenberg \cite{Izzo-Levenberg:2019}, we will establish results on pairs of polynomially convex sets whose union has nontrivial hull without analytics discs (Theorems~\ref{no-discs3} and~\ref{no-discs4} and Corollary~\ref{no-discs-gen}).

All of the examples we will produce rely on the same fundamental idea originally introduce by Wermer \cite{Wermer:1955} to construct an arc with nontrivial hull in $\C^3$ and developed further by Walter Rudin \cite{Rudin:1956} to construct a Cantor set with nontrivial hull, and an arc with nontrivial hull, in $\C^2$.
In the next section, in addition to defining terminology and notation (some of which has already been used above), we will introduce a certain algebra $A_J$ on which these constructions depend and summarize the fact we will need about $A_J$.

Our main results concerning pairs of polynomially convex sets whose union is not rationally convex are treated in Section~\ref{nonrat-unions}.  The approximation of simple closed curves by locally polynomially convex, nonrationally convex simple closed curves is treated in Section~\ref{path-curves}.  The results concerning the union of the closed unit ball and a polynomially convex arc are in Section~\ref{ball-n-chain}.  The hulls without analytic discs are in the concluding Section~\ref{no-discs}.

It is with a mixture of joy and sorrow that I dedicate this paper to the memory of Peter Duren.  Sorrow, of course, that he is no longer with us; joy that I had the privilege of knowing him.  

%
%

\section{Preliminaries}\label{prelim}

By definition an \emph{arc} is a space homeomorphic to the closed unit interval, and a \emph{simple closed curve} is a space homeomorphic to the unit circle.  For convenience  we will also use the term \emph{arc} (or \emph{simple closed curve}) to refer to a topological embedding whose domain is the closed unit interval (or the unit circle).  By a 
\emph{Cantor set} we mean any space that is homeomorphic to the usual mid-thirds Cantor set.  We will make repeated use of the well known characterization of Cantor sets as the compact, totally disconnected, metrizable spaces without isolated points.

We call a space \emph{simply coconnected} if its first \v Cech cohomology group with integer coefficients is zero.  It is well known that a compact set in the plane is simply coconnected if and only if its complement is connected.

For
$X$ a compact Hausdorff space, we denote by $C(X)$ the algebra of all continuous complex-valued functions on $X$ with the supremum norm
$ \|f\|_{X} = \sup\{ |f(x)| : x \in X \}$.  A \emph{uniform algebra} on $X$ is a closed subalgebra of $C(X)$ that contains the constant functions and separates
the points of $X$.  The maximal ideal space of a uniform algebra $A$ will be denoted by $\ma$.

For a compact set $X$ in $\C^N$, the \emph{polynomial hull} $\h X$ of $X$ is defined by
$$\h X=\{z\in\C^N:|p(z)|\leq \max_{x\in X}|p(x)|\
\mbox{\rm{for\ all\ polynomials}}\ p\},$$
and the
\emph{rational hull} $\hr X$ of $X$ is defined by
$$\hr X = \{z\in\C^N: p(z)\in p(X)\ 
\mbox{\rm{for\ all\ polynomials}}\ p
\}.$$
An equivalent formulation of the definition of $\hr X$ is that $\hr X$ consists precisely of those points $z\in \C^N$ such that every polynomial that vanishes at $z$ also has a zero on $X$.
The set $X$ is said to be \emph{polynomially convex} if $\h X=X$ and \emph{rationally convex} if $\hr X=X$.  We say that a polynomial hull $\h X$ (or rational hull $\hr X$) is \emph{nontrivial} if $\hh X$ (or $\hhr X$) is nonempty.

We denote by 
$P(X)$ the uniform closure on $X\subset\C^N$ of the polynomials in the complex coordinate functions $z_1,\ldots, z_N$, and we denote by $R(X)$ the uniform closure of the rational functions  holomorphic on (a neighborhood of) $X$. 
Both $P(X)$ and $R(X)$ are uniform algebras, and
it is well known that the maximal ideal space of $P(X)$ can be naturally identified with $\h X$, and the maximal ideal space of $R(X)$ can be naturally identified with $\hr X$.

As usual, we denote the Gelfand transform of a function $f$ by~$\hat f$.
Given a uniform algebra $A$ and a compact subset $E$ of the maximal ideal space $\ma$ of $A$, the \emph{$A$-convex hull} of $E$ is the set
$$\h E=\{\phi \in \ma: |\hat f(\phi)| \leq \|f\|_E\ {\rm for\ all\ } f\in A \}.$$
The $A$-convex hull of a compact set can be thought of as an abstract polynomial hull.  Indeed, when $E\subset X$ are compact sets in $\C^N$, and $A$ is taken to be $P(X)$, the $A$-convex hull of $E$ coincides with the polynomial hull of $E$.  There is also an analogous abstract rational hull defined as follows.

\bdefn
Given a uniform algebra $A$ and a compact subset $E$ of the maximal ideal space $\ma$ of $A$, we define the {\em $A$-rational hull $\hr E$ of $E$} to be the set
$$\hr E=\{\phi \in \ma: \hat f(\phi)\in \hat f(E)\ {\rm for\ all\ } f\in A \}.$$
We say that $E$ is {\em $A$-rationally convex} if $\hr E=E$.
\edefn

An equivalent formulation of the definition of the $A$-rational hull $\hr E$ is that $\hr E$ consists precisely of those points $\phi \in \ma$ such that if $f\in A$ satisfies $\hat f(\phi)=0$, then $\hat f$ has a zero on $E$.  When $E\subset X$ are compact sets in $\C^N$, and $A$ is taken to be $P(X)$, the $A$-rational hull of $E$ coincides with the usual rational hull of $E$.

For $A$ a uniform algebra on a compact space $X$, the \emph{Shilov boundary} for $A$ is the smallest closed subset $\Gamma$ of $X$ such that $\| f \|_{\Gamma}=\|f\|_X$.  Thus, the Shilov boundary for $A$ is the smallest closed subset whose $A$-convex hull is $\ma$.  (For the existence of the Shilov boundary see \cite[pp.~94--95]{Browder:1969}, \cite[p.~9]{Gamelin:1984}, or \cite[p.~15]{Stout:2007}.)

By an \emph{analytic disc} in $\C^N$, we mean an injective holomorphic map from an open disc in the complex plane into $\C^N$.
By the statement that a subset $S$ of $\C^N$ contains no analytic discs, we mean that there is no analytic disc in $\C^N$ whose image is contained in $S$.

Given a compact planar set $J$, we denote by $A_J$ the algebra of continuous complex-valued functions on the Riemann sphere $S^2$ that are holomorphic on $S^2\sm J$.  We summarize here the facts we need about $A_J$.  More detail can be found in many texts, for instance \cite[pp.~53--54]{Stout:2007}.  In general, $A_J$ need not contain nonconstant functions.  However, whenever $A_J$ does contain a nonconstant function, then there are three functions in $A_J$ that separate points on $S^2$, and consequently, $A_J$ is a uniform algebra.  Indeed, if $g$ is a nonconstant function in $A_J$ and $a_1$ and $a_2$ are points in $\C\sm J$ such that $g(a_1)\neq g(a_2)$, then, regarding $\big((g-g(a_j)\bigl)/(z-a_j)$, $j=1,2$ as functions in $A_J$, it is well known (and easily verified) that the three functions $g, \big((g-g(a_1)\bigl)/(z-a_1), \big((g-g(a_2)\bigl)/(z-a_2)$ separate points on $S^2$.  The algebra $A_J$ contains a nonconstant function whenever $J$ has positive planar measure, for then the function $g$ defined by
$$g(\zeta)=\int\!\!\!\int_J\, \frac{dx\,dy}{z-\zeta}$$
is such a function.  It is immediate from the maximum principle that every function in $A_J$ takes its maximum modulus on $J$, so the $A_J$-convex hull of $J$ is $S^2$.  In fact, for every function $f\in A_J$, the inclusion $f(S^2)\subset f(J)$ holds, so the $A_J$-rational hull of $J$ is $S^2$.  For a proof, see \cite[p.~54]{Stout:2007}.  It follows immediately that given $\row fN\in A_J$ and letting $S^2\rightarrow\C^N$ be defined by $\pi(z)=\bigl(f_1(z),\ldots,f_N(z)\bigr)$, the inclusion $\pi(S^2)\subset \hr{\pi(J)}$ holds, and hence $\pi(J)$ is not rationally convex whenever $\pi(S^2)\not\subset \pi(J)$.

\bprop\label{Shilov1}
Let $J\subset \C$ be a compact set such that $A_J$ contains a nonconstant function, and let $\Gamma$ denote the Shilov boundary for $A_J$.  Then $A_\Gamma=A_J$.
\eprop

\bpf
This is \cite[exercise II.1(a)]{Gamelin:1984}.  
\epf

\bprop\label{measure}
Let $J\subset \C$ be a compact set such that $A_J$ contains a nonconstant function, and let $\Gamma$ denote the Shilov boundary for $A_J$.  Then 
$J\sm\Gamma$ has planar measure zero.
\eprop

\bpf
Let $x_0\in J$ be arbitrary, and suppose that the intersection of each disc centered at $x_0$ with $J$ has positive measure.  We will show that then $x_0$ is in the Shilov boundary for $A_J$.  The lemma follows.

Let $r>0$ be arbitrary.  Choose a compact subset $E$ of $J\cap D(x_0,r)$ of positive measure.  Define a function on $S^2$ by
$$f(\zeta)=\int\!\!\! \int_E \frac{dx\,dy}{z-\zeta}.$$
It is well-known that $f$ is in $A_E\subset A_J$ and $f$ is nonconstant in the unbounded component of the complement of $E$.  Thus $|f(\zeta)|<\|f\|_{S^2}$ for every point $\zeta\in S^2\sm D(x_0,r)$.  Consequently, the Shilov boundary for $A_J$ must intersect $D(x_0,r)$.  Since $r>0$ was arbitrary, this shows that $x_0$ is in the Shilov boundary for $A_J$.
\epf

\bcor\label{Shilov2}
Let $J\subset \C$ be a compact set that is locally of positive planar measure.  Then $J$ is the Shilov boundary for $A_J$.
\ecor

%
%

\section{Nonrationally convex unions}\label{nonrat-unions}

In this section we establish several results regarding pairs of polynomially convex sets whose union is not even rationally convex.  We will begin with a result giving, in the abstract uniform algebra setting, minimal sets with nontrivial hull.  
Note that the condition that a nonrationally convex set $\Gamma$ is a minimal nonpolynomially convex set is stronger than the condition that $\Gamma$ is the union of two polynomially convex sets.
We will then give a general result about pairs of polynomially convex sets in $\C^3$ whose union is not rationally convex from which we will obtain several concrete examples.  After that we will give some examples in $\C^2$.  We will conclude with a general result about embedding a pair of spaces in some $\C^N$ so as to be polynomially convex but have nonrationally convex union.

By the complement and interior of a planar set, we mean its complement and interior relative to the Riemann sphere.

\bthm\label{abstract-case}
Let $J$ be a compact planar set such that $A_J$ contains a nonconstant function and such that $J$ is the boundary of each component of the complement of $J$.  Let $\Gamma$ denote the Shilov boundary for $A_J$.  Then the $A_\Gamma$-rational hull of $\Gamma$ is $S^2$, but every proper closed subset $K$ of $\Gamma$ is $A_\Gamma$-convex and satisfies $\ol{A_\Gamma|K}=C(K)$.
\ethm

Note the analogy with the disc algebra: the polynomial hull of the unit circle in the plane is the closed unit disc 
while every proper closed subset of the circle is polynomially convex and the polynomials are dense in the continuous functions there.  A similar analogy holds with the big disc algebra.  However, in contrast to the situation with the disc algebra and big disc algebra, $\Gamma$ has $A_\Gamma$-{\em rational}-hull, not just $A_\Gamma$-polynomial hull.

There are examples of the situation in Theorem~\ref{abstract-case} with $\Gamma$ a Cantor set or an arc.  To see this, simply take $J$ to be a Cantor set or arc that is locally of positive planar measure.  Then $\Gamma=J$ by Corollary~\ref{Shilov2}. 

For the proof of Theorem~\ref{abstract-case} we will need the following topological lemma.

\blem\label{connected-complement}
Let $J$ be a compact planar set such that $J$ is the boundary of each component of the complement of $J$.  Let $K$ be a proper closed subset of $J$.  Then $K$ has empty interior and connected complement.
\elem

\bpf
Obviously $J$ has empty interior, and hence the same is true of $K$.  Since $J$ has empty interior, each component of the complement of $K$ intersects the complement of $J$ and hence contains a component of the complement of $J$.  Let $x$ be a point of $J\sm K$.  Then some disc $\Delta$ centered at $x$ is contained in $\C\sm K$.  This disc intersects every component of the complement of $J$ and hence intersects every component of the complement of $K$.  Consequently, the disc $\Delta$ is contained in every component of the complement of $K$.  Thus the complement of $K$ consists of a single component.
\epf

\bpf[Proof of Theorem~\ref{abstract-case}]
By Proposition~\ref{Shilov1}, $A_\Gamma=A_J$, and as mentioned just before the proposition, for any compact planar set $E$ for which $A_E$ contains a nonconstant function, 
the $A_E$-rational hull of $E$ is $S^2$.

Now let $K$ be a proper closed subset of $\Gamma$.  Since $\Gamma$ is the Shilov boundary for $A_\Gamma$, there is a function $f\in A_\Gamma$ such that $\|f\|_K< \|f\|_{S^2}=\| f\|_\Gamma$.  Choose a point $z_0\in \Gamma\sm K$ such that $|f(z_0)|=\| f\|_{S^2}$.  Set $\vep=|f(z_0)|-\|f\|_K>0$.  By \cite[Theorem~II.1.8]{Gamelin:1984}, there is a function $g\in A_\Gamma$ such that $g$ is holomorphic in a neighborhood of $z_0$ and $\|f-g\|_{S^2}< \vep/2$.  Then $|g(z_0)|> \|g\|_K$.  In particular, $g$ never takes the value $g(z_0)$ on the $A_\Gamma$-convex hull $\h K$ of $K$.  Therefore, the function $\bigl(g-g(z_0)\bigr)|K$ is invertible as an element of the uniform algebra $\ol{A_\Gamma|K}$.

Since $g-g(z_0)$ is holomorphic in a neighborhood of $z_0$ and vanishes at $z_0$, we can regard $\bigl(g-g(z_0)\bigr)/(z-z_0)$ as a continuous function on $S^2$ that is holomorphic on $S^2\sm \Gamma$ and thus belongs to $A_\Gamma$.  On $K$ we have
$$\left( \frac{g-g(z_0)}{z-z_0} \right) \bigl(g-g(z_0)\bigr)^{-1} = \frac{1}{z-z_0}$$
so $(z-z_0)^{-1}$ is in $\ol{A_\Gamma|K}$.  Since by Lemma~\ref{connected-complement}, $K$ has empty interior and connected complement, Lavrentiev's theorem \cite[Theorem~II.8.7]{Gamelin:1984} gives that every continuous function on $K$ can be approximated uniformly by polynomials.  By Runge's theorem, every polynomial can be approximated uniformly on $K$ by rational functions whose only pole is at $z_0$.  Thus 
$\ol{A_\Gamma|K}=C(K)$ and, in particular, $K$ is $A_\Gamma$-convex.
\epf

We turn now to hulls in $\C^3$.

\bthm\label{union-of-poly-convex}
Let $J$ be a compact planar set with empty interior such that $A_J$ contains a nonconstant function.  Let $K_1$ and $K_2$ be compact sets whose union is $J$ neither of which contains the Shilov boundary for $A_J$ and each of which has connected complement.  Then there is an embedding $\pi$ of $J$ into $\C^3$ such that $\pi(J)$ has nontrivial rational hull, but $\pi(K_1)$ and $\pi(K_2)$ are each polynomially convex and satisfy $P(\pi(K_j))=C(\pi(K_j))$, $j=1,2$.
\ethm

This theorem can be reformulated in the following more intrinsic manner.

\bcor\label{intrinsic-form}
Let $(J,K_1,K_2)$ be a triple with $J$ an uncountable, compact space of topological dimension at most 1 that embeds in the plane and $K_1$ and $K_2$ closed subspaces of $J$ whose union is $J$ and such that each of $K_1$ and $K_2$ is simply coconnected and each of $J\sm K_1$ and $J\sm K_2$ is uncountable.  Then there is an embedding $\pi$ of $J$ 
into $\C^3$ such that $\pi(J)$ has nontrivial rational hull, but $\pi(K_1)$ and $\pi(K_2)$ are each polynomially convex and $P(\pi(K_j))=C(\pi(K_j))$, $j=1,2$. 
 \ecor

Before proving Theorem~\ref{union-of-poly-convex} and Corollary~\ref{intrinsic-form}, we present several other corollaries.  In all these corollaries the polynomially convex sets can be chosen to be sets on which the polynomials are dense in the continuous functions.  This has been omitted from the statements to avoid unnecessary repetition.

\bcor\label{disjoint-Cantor-sets}
There exists a Cantor set in $\C^3$ that has nontrivial rational hull and is the union of two disjoint polynomially convex Cantor sets.
\ecor

\bpf
Take $J$ to be a Cantor set that is locally of positive planar measure, take $K_1$ and $K_2$ to be a separation of $J$, and apply Theorem~\ref{union-of-poly-convex}.
\epf

Taking $J$ in Theorem~\ref{union-of-poly-convex} to be an arc that is locally of positive planar measure we obtain:

\bcor\label{arc-union}
There exists an arc in $\C^3$ that has nontrivial rational hull and is the union of two  polynomially convex arcs.
\ecor

Note that the subarcs can be chosen to have only a point in common, or alternatively, they can be chosen so as to overlap so that the nonrationally convex arc is locally polynomially convex.

\bcor\label{disjoint-arcs-C3}
There exist two disjoint polynomially convex arcs in $\C^3$ whose union is not rationally convex.
\ecor

\bcor\label{scc-not-rat-convex}
There exists a simple closed curve in $\C^3$ that has nontrivial rational hull and is the union of two polynomially convex arcs.
\ecor

Of course the existence of two polynomially convex arcs whose union is a nonpolynomially convex simple closed curve is trivial; just take two arcs in the plane whose union is the unit circle.  However, that the union can fail to be rationally convex as in the above corollary does not seem to be obvious.
Note that the arcs in the corollary can be chosen so that their interiors cover the nonrationally convex simple closed curve and hence the nonrationally convex simple closed curve is locally polynomially convex.

In preparation for the proof of Theorem~\ref{union-of-poly-convex} we establish a lemma.  (For $A$ a uniform algebra on a compact space $X$, a point $x_0$ of $X$ is said to be a \emph{peak point} for $A$ if there is a function in $A$ that \emph{peaks} at $x_0$, that is a function $f$ in $A$ such that $f(x_0)=1$ and $|f(x)|<1$ for every $x\in X\sm \{x_0\}$.)

\blem\label{special-function}
Let $A$ be a uniform algebra on a compact space $X$.  Let $x_0\in X$ be a peak point for $A$, and let $K$ be a compact subset of $X$ that does not contain $x_0$.  Then there is a function $f$ in $A$ such that 
\begin{enumerate}
\item[(i)] $f(x_0)=1$
\item[(ii)] $\Re f\geq 0$ everywhere on $X$
\item[(iii)]  $\Re f \leq 1/4$ on $K$.
\end{enumerate}
\elem

\bpf
Choose a function $h\in A$ that peaks at $x_0$.  Let $r=\|h\|_{K_1}<1$, and let $\alpha$ be the real number such that the conformal automorphism of the disc defined by $\varphi_\alpha(z)=({z-\alpha})/({1-\alpha z})$ takes $r$ to $-1/2$.  Then the function $f=\bigl(1+(\varphi_\alpha\circ h)\bigr)/2$ has the desired properties.
\epf

\bpf[Proof of Theorem~\ref{union-of-poly-convex}]
Since neither of $K_1$ and $K_2$ contains the Shilov boundary for $A_J$, there exist peak points $a_1$ and $a_2$ for $A_J$ with $a_1\notin K_1$ and $a_2\notin K_2$.  By Lemma~\ref{special-function} there are functions $f_1$ and $f_2$ in $A_J$ such that for each $j=1,2$, we have
\begin{enumerate}
\item[(i)] $f_j(a_j)=1$
\item[(ii)] $\Re f_j\geq 0$ everywhere on $J$
\item[(iii)]  $\Re f_j \leq 1/4$ on $K_j$.
\end{enumerate}
Set $f=f_1-f_2$.  Then
\begin{enumerate}
\item[(i)] $\Re f(a_1) \geq 3/4 $
\item[(ii)] $\Re f\leq 1/4$ on $K_1$
\item[(iii)]  $\Re f(a_2) \leq \, -3/4$
\item[(iv)]  $\Re f \geq \, -1/4$ on $K_2$.
\end{enumerate}
Repeated application of \cite[Theorem~II.1.8]{Gamelin:1984} yields a function $g\in A_J$ such that $g$ is holomorphic in neighborhoods of each of $a_1$ and $a_2$ and satisfies $\|f-g\|_{S^2} < 1/4$.  Then
\begin{enumerate}
\item[(i)] $\Re g(a_1) > 1/2 $
\item[(ii)] $\Re g< 1/2$ on $K_1$
\item[(iii)]  $\Re g(a_2) < \, -1/2$
\item[(iv)]  $\Re g > \, -1/2$ on $K_2$.
\end{enumerate}

Note that $\bigl(g-g(a_1)\bigr) / (z-a_1)$ and $\bigl(g-g(a_2)\bigr) / (z-a_2)$ can be regarded as functions in $A_J$.  Let $\pi:S^2\to \C^3$ be defined by
$\pi(z)= \bigl(g(z), (g(z)-g(a_1)\bigr) / (z-a_1), (g(z)-g(a_2)\bigr) / (z-a_2) \bigr)$.
Then (see the end if Section~\ref{prelim}) $\pi$ is injective and hence is an embedding,  and  $\pi(S^2$) is contained in the rational hull of $\pi(J)$, so the rational hull of $\pi(J)$ is nontrivial.

Let $A$ denote the uniform algebra on $S^2$ generated by the three functions $g$, $\bigl(g-g(a_1)\bigr) / (z-a_1)$, and $\bigl(g-g(a_2)\bigr) / (z-a_2)$.
Because on $K_1$ we have $\Re g< 1/2< g(a_1)$, the function 
$1/\bigl(z-g(a_1)\bigr)$ can be approximated uniformly on $g(K_1)$ by polynomials.  Consequently, the function $\bigl(g-g(a_1)\bigr)|K_1$ is invertible as an element of the uniform algebra $\ol{A|K_1}$.  (Alternatively, one can obtain this by noting that because on $K_1$ we have $\Re g< 1/2< g(a_1)$, the function $g$ never takes the value $g(a_1)$ on the $A$-convex hull of $K_1$.)  On $K_1$ we have
$$\left( \frac{g-g(a_1)}{z-a_1} \right) \bigl(g-g(a_1)\bigr)^{-1} = \frac{1}{z-a_1}$$
so $(z-a_1)^{-1}$ is in $\ol{A|K_1}$.  Now exactly as in the proof of Theorem~\ref{abstract-case}, we get that $\ol{A|K_1}=C(K_1)$.  It follows that $P(\pi(K_1))=C(\pi(K_1))$ and hence $\pi(K_1)$ is polynomially convex.  The proof that the same holds for $K_2$ is similar.
\epf


The proof of Corollary~\ref{intrinsic-form} uses the following lemma which I learned from Stout.

\blem
Given Cantor sets $E$ and $E^*$ in $\C$, there exists a homeomorphism $h:\C\to\C$ such that $h(E)=E^*$.
\elem

\bpf
Every Cantor set in the plane is contained in a simple closed curve in the plane.  Choose simple closed curves $J$ and $J^*$ that contain $E$ and $E^*$, respectively.  Given any two Cantor sets $G$ and $G^*$ in the real line $\R$, there is a homeomorphism of $\R$ onto itself that maps $G$ onto $G^*$.  It follows that there is a homeomorphism $g:J\rightarrow J^*$ such that $g(E)=E^*$.  By the Schoenflies theorem \cite[Ch.~5, Corollary~9.25]{Hall-Spencer:1955}, $g$ extends to a homeomorphism $h:\C\to\C$.
\epf

\bpf[Proof of Corollary~\ref{intrinsic-form}]
By hypothesis each of $J\sm K_1$ and $J\sm K_2$ is an uncountable open set of 
$J$.  Every open set in a metrizable space is a countable union of closed subsets of the space.  Consequently, each of $J\sm K_1$ and $J\sm K_2$ contains an uncountable closed, and hence compact, subset of $J$.  Every uncountable compact metrizable space contains a Cantor set.  Thus we can choose Cantor sets $G_1$ and $G_2$ in $J\sm K_1$ and $J\sm K_2$, respectively.  Note that then $G_1\cup G_2$ is also a Cantor set.

Choose an embedding $g:J \to \C$.  Choose a Cantor set $L$ in $\C$ that is locally of positive planar measure.  By the preceding lemma, there is a homeomorphism $h:\C\to \C$ such that $h\bigl(g(G_1\cup G_2)\bigr)=L$.  Let $f:J\to \C$ be the embedding given by $f=h\circ g$.  Then each of $f(G_1)$ and $f(G_2)$ is a Cantor set of positive planar measure.  Consequently, each of $f(J)\sm f(K_1)$ and $f(J)\sm f(K_2)$ has positive planar measure and hence intersects the Shilov boundary for $A_{f(J)}$, by
Theorem~\ref{measure}.

Since $J$ has topological dimension at most 1, and hence the same is true of $K_1$ and $K_2$, each of $f(K_1)$ and $f(K_2)$ has empty interior in $\C$.  The simply coconnectivity of $K_1$ and $K_2$ implies that each of $f(K_1)$ and $f(K_2)$ has connected complement in $\C$.

Now Theorem~\ref{union-of-poly-convex} applies to yield the desired conclusion.
\epf

We turn next to examples in $\C^2$.  

\bthm\label{two-Cantor-sets}
There exists a Cantor set in $\C^2$ that has nontrivial rational hull and is the union of two polynomially convex Cantor sets.
\ethm

Whether the two polynomially convex Cantor sets in this theorem can be taken to be disjoint, as was shown in the case of sets in $\C^3$, remains open.

\bcor\label{C2-two-arcs}
There exist two polynomially convex arcs in $\C^2$ whose union is not rationally convex.
\ecor

\bpf
By the preceding theorem there exist polynomially convex Cantor sets $K_1$ and $K_2$ in $\C^2$ whose union $E$ is not rationally convex.  Every Cantor set in a Euclidean space is contained in an arc \cite[Theorem~p.~57]{Whyburn:1932}.  Therefore, by \cite[Corollary~1.2]{Izzo:2022} there are polynomially convex arcs $J_1$ and $J_2$ containing $K_1$ and $K_2$, respectively such that each of $J_1\sm K_1$ and $J_2\sm K_2$ is a countable union of smooth arcs.  Since the rational hull $\hr {J_1\cup J_2}$ of $J_1\cup J_2$ contains $\hhr E$, and the later set contains a smooth 2-manifold, it follows that $J_1\cup J_2$ is not rationally convex.
\epf

Whether Corollary~\ref{C2-two-arcs} can be strengthened to show that Corollaries~\ref{arc-union}--\ref{scc-not-rat-convex} hold with $\C^2$ in place of $\C^3$ remains open.

\bpf[Proof of Theorem~\ref{two-Cantor-sets}]
Let $K_1$ and $K_2$ be two Cantor sets in $\C$ that are locally of positive planar measure.  Let $J=K_1\cup K_2$. 
Define functions $f_1$ and $f_2$ on $S^2$ by 
$$f_1(\zeta)=\int\!\!\! \int_{K_1} \frac{dx\,dy}{z-\zeta}\qquad {\rm and}\qquad f_2(\zeta)=\int\!\!\! \int_{K_2} \frac{dx\,dy}{z-\zeta}.$$
As noted in Section~\ref{prelim}, each $f_j$ is in $A_{K_j}$.  Furthermore, $f'_j(\infty)=\lim_{j\rightarrow\infty} \zeta f(\zeta)=\int\!\!\! \int_{K_j}  dx\, dy\neq 0$.   Thus by replacing $K_2$ by a translate of itself that is sufficiently far away from $K_1$, we can arrange to have $f_1$ be one-to-one on $K_2$ and $f_2$ be one-to-one on $K_1$.  Note that of course $f_1$ and $f_2$ belong to $A_J$.  

Define $\pi:J\to \C^2$ by $\pi(z)=\bigl(f_1(z), f_2(z)\bigr)$.  
Note that $\pi$ maps each $K_j$ homeomorphically onto $\pi(K_j)$.  Consequently, each of $\pi(K_1)$, 
$\pi(K_2)$, and $\pi(J)$ is a Cantor set.
Furthermore, because the coordinate function $z_2$ is one-to-one on $K_1$, Lavrentiev's theorem gives that $P(\pi(K_1))=C(\pi(K_1))$.  The same holds with the roles of 1 and 2 reversed.
In particular, each of $\pi(K_1)$ and $\pi(K_2)$ is polynomially convex.
As noted in Section~\ref{prelim}, $\pi(S^2)$ is contained in the rational hull of $\pi(J)$, and hence the rational hull of $\pi(J)$ is nontrivial.
\epf

We conclude this section by showing that
it follows from the existence of disjoint polynomially convex Cantor sets whose union has nontrivial rational hull that any uncountable compact subspace of a Euclidean space can be embedded in some $\C^N$ as a nonrationally convex set that is the union of two polynomially convex sets.  First we consider an abstract form of this result that applies to more general compact spaces.

\bthm\label{general-unions}
Let $(J,K_1,K_2)$ be a triple with $J$ a compact Hausdorff space and $K_1$ and $K_2$ closed subspaces of $J$ whose union is $J$ and such that each of $J\sm K_1$ and $J\sm K_2$ contains a Cantor set.  Then there exists a uniform algebra $A$ on $J$ such that $J$ has nontrivial $A$-rational hull, but
each of $K_1$ and $K_2$ is $A$-convex and $\ol{A|K_j}=C(K_j)$, $j=1,2$.
\ethm
 
 \bpf
Choose Cantor sets $G_1$ and $G_2$ in $J\sm K_1$ and $J\sm K_2$, respectively.  Then Corollary~\ref{disjoint-Cantor-sets} yields a uniform algebra $B$ on $G_1\cup G_2$ such that $G_1\cup G_2$ has nontrivial $B$-rational hull, but $G_1$ and $G_2$ are each $B$-convex.  Furthermore, $B$ can be chosen so that $\ol{B|G_1}=C(G_1)$ and $\ol{B|G_2}=C(G_2)$.  Let $A$ be the uniform algebra on $J$ defined by $A=\{f\in C(J): f|(G_1\cup G_2)\in B\}$.
 
 Each of $J$ and $\mb$ can be regarded as subsets of $\ma$ in standard ways.  Then $\mb$ is the $A$-convex hull of $G_1\cup G_2$ in $\ma$, and it follows that $J\cap \mb=G_1\cup G_2$.  Since the $B$-rational hull of $G_1\cup G_2$ is easily seen to be contained in the $A$-rational hull of $G_1\cup G_2$, we obtain that the $A$-rational hull of $G_1\cup G_2$ is not contained in $J$.  Consequently, $J$ has nontrivial $A$-rational hull.
 
Since $\ol{B|G_2}=C(G_2)$, it follows from the Tietze extension theorem that the functions in $C(K_1)$ whose restrictions to $G_2$ lie in $B|G_2$ form a dense subset of $C(K_1)$.  Thus to show that $\ol{A|K_1}=C(K_1)$, it suffices to show that each function $f\in C(K_1)$ that satisfies $f|G_2\in B|G_2$ extends to a member of $A$.  To do so, choose $g\in B\subset C(G_1\cup G_2)$ such that $g|G_2=f$, note that the formula 
$$h(x)= \begin{cases}
f(x), & x\in K_1\\
g(x), & x\in G_1
\end{cases}$$
yields a well-defined continuous function on $K_1\cup G_1$, extend $h$ to a continuous function $\tilde h$ on $J$, and note that then $\tilde h$ is in $A$ and $\tilde h|K_1=f$.  Thus $\ol{A|K_1}=C(K_1)$.  The verification that $\ol{A|K_2}=C(K_2)$ is similar.  Finally, the equation $\ol{A|K_j}=C(K_j)$ implies $A$-convexity of $K_j$.
 \epf
 
 \bcor\label{n-dimensional-case}
 Let $(J,K_1,K_2)$ be a triple with $J$ a compact subspace of $\R^N$ and $K_1$ and $K_2$ closed subspaces of $J$ whose union is $J$ and such that each of $J\sm K_1$ and $J\sm K_2$ is uncountable.  Then there exists an embedding $\pi$ of $J$ into $\C^{N+4}$ such that $\pi(J)$ has nontrivial rational hull, but
each of $\pi(K_1)$ and $\pi(K_2)$ is polynomially convex and $P(\pi(K_j))=C(\pi(K_j))$, $j=1,2$.
 \ecor
 
Since every compact metrizable space of topological dimension $m$ embeds in $\R^{2m+1}$ \cite[Theorem~V.2]{Hurewicz-Wallman:1948}, Theorem~\ref{n-dimensional-case} show, in particular, that there is a homeomorphic copy of every uncountable such space in $\C^{2m+5}$ that has nontrivial rational hull but is the union of two polynomially convex sets.  As another immediate consequence of Corollary~\ref{n-dimensional-case}, if $K_1$ and $K_2$ are any two uncountable compact subsets of $\R^N$, then there are disjoint polynomially convex sets in $\C^{N+4}$ that are homeomorphic to $K_1$ and $K_2$ but whose union is not rationally convex.

\bpf
As in the proof of Corollary~\ref{intrinsic-form}, we can choose Cantor sets $G_1$ and $G_2$ in $J\sm K_1$ and $J\sm K_2$, respectively.  As in the proof of Theorem~\ref{general-unions}, let $B$ be a uniform algebra on $G_1\cup G_2$ such that $G_1\cup G_2$ has nontrivial $B$-rational hull, but $G_1$ and $G_2$ are each $B$-convex, and $\ol{B|G_j}=C(G_j)$, $j=1,2$.    By Corollary~\ref{disjoint-Cantor-sets}, the uniform algebra $B$ can be chosen so as to be generated by three functions $f_1,f_2,f_3$.  Set 
$A=\{f\in C(J): f|(G_1\cup G_2)\in B\}$.

Extend each of $f_1,f_2,f_3$ to continuous complex-valued functions $\tilde f_1,\tilde f_2,\tilde f_3$ on $J$.  Let $x_1,\ldots, x_N$ denote the real coordinate functions of $\R^N$.  Choose a continuous real-valued function $\rho$ on $J$ whose zero set is precisely $G_1\cup G_2$.  Then the $N+4$ functions $\tilde f_1,\tilde f_2,\tilde f_3, \rho, \rho x_1,\ldots \rho x_N$ generate the uniform algebra $A$ by \cite[Lemma~3.8]{ISW:2016}.  The proof of Theorem~\ref{general-unions} shows that $J$ has nontrivial $A$-rational hull and $\ol{A|K_j}=C(K_j)$, $j=1,2$.  It follow that the map $\pi:J\to \C^{N+4}$ whose component functions are the functions $\tilde f_1,\tilde f_2,\tilde f_3, \rho, \rho x_1,\ldots \rho x_N$ has the properties asserted in the theorem.
\epf

%
%

\section{Approximation by pathological curves}\label{path-curves}

In this section we prove that every simple closed curve in $\C^N$, $N\geq 3$, can be approximated uniformly
by a nonrationally convex simple closed curve that is locally polynomially convex, and is, in fact, covered by the interiors of \emph{two} polynomially convex arcs.  We will use $\|\cdot\|_\infty$ to denote the supremum norm.

\bthm\label{path-curve}
Let $\gamma$ be a simple closed curve in $\C^N$, $N\geq 3$.  Given $\vep>0$, there exists a simple closed curve $\sigma$ satisfying $\|\gamma -\sigma\|_\infty<\vep$ that is not rationally convex but is the union of two polynomially convex arcs $\sigma_1$ and $\sigma_2$ satisfying $P(\sigma_j)=C(\sigma_j)$, $j=1,2$, and whose interiors cover $\sigma$.
\ethm

Since every simple closed curve in Euclidean space can be approximated uniformly by a smooth, and hence rectifiable, simple closed curve, the above theorem is a consequence of the following sharper result for rectifiable simple closed curves.

\bthm\label{rectifiable-approximated}
Given a rectifiable simple closed curve $\gamma$ in $\C^N$, $N\geq 3$, given $\vep>0$, and given an open ball $B$ of $\C^N$ that intersects $\gamma$, there is a simple closed curve $\gamma_a$ that satifies $\|\gamma-\gamma_a\|_\infty<\vep$ and $\gamma_a\sm B=\gamma \sm B$ such that $\gamma_a$ is not rationally convex but there are two arcs $\gamma_1$ and $\gamma_2$ that are polynomially convex and satisfy $P(\gamma_j)=C(\gamma_j)$ and whose interiors cover $\gamma_a$.
\ethm

In proving that certain arcs or simple closed curves are not polynomially or rationally convex, the following lemma is useful.

\blem\label{convexity-of-subsets}
If an arc or simple closed curve $\gamma$ is polynomially convex, then so is every subset of $\gamma$.  The same statement holds with polynomial convexity replaced by rational convexity.
\elem

The lemma is a consequence of the following theorem, a slight generalization of \cite[Corollary~III.8.3]{Gamelin:1984}.  

\bthm\label{gen-Rossi-consequence}
Let $A$ be a uniform algebra with 
maximal ideal space $\ma$ and Shilov boundary $\Gamma_A$.  Suppose 
$\ma\sm\Gamma_A$ is homeomorphic to a subset of a one-dimensional manifold. Then $\ma=\Gamma_A$.
\ethm

\bpf
The proof is similar to the proof of \cite[Corollary~III.8.3]{Gamelin:1984}.  Assume to get a contradiction that $\ma\sm\Gamma_A$ is nonempty, and fix a point $x$ in $\ma\sm\Gamma_A$.  Let $U$ be a neighborhood of $x$ contained in $\ma\sm\Gamma_A$.  By Rossi's local maximum modulus principle \cite[Theorem~III.8.2]{Gamelin:1984}, if $f$ is a function in $A$ whose Gelfand transform $\hat f$ vanishes identically on $\partial U$, then $\hat f$ must vanish also at $x$.  Since $U$ can be chosen so that its boundary consists of at most two points, this is a contraction.
\epf

\blem
If the polynomial hull of a compact set $E\subset \C^N$ is contained in an arc or simple closed curve, then $E$ is polynomially convex.  The same statement holds with polynomial convexity replaced by rational convexity.
\elem

\bpf
Apply the theorem to the uniform algebra $P(E)$ or $R(E)$.
\epf

Lemma~\ref{convexity-of-subsets} is an immediate consequence of the last lemma.

For the reader's convenience we quote here two other results we will use.

\begin{theorem}\cite[Theorem~1.6]{Izzo-Stout:2022b}\label{from-Stout-paper}
Let $Y$ be a compact polynomially convex subset of $\C^N$, and let $\Gamma$ be a subset of $\C^N$ such that $Y\cup \Gamma$ is compact and such that for every neighborhood $U$ of $Y$ in $\C^N$, the set $\Gamma\sm U$ is contained in a compact connected set of finite length.  Suppose also that the map $\check H^1(Y\cup \Gamma; \Z)\rightarrow \check H^1(Y;\Z)$ induced by the inclusion $Y\hookrightarrow Y\cup \Gamma$ is a monomorphism.  Then $Y\cup \Gamma$ is  polynomially convex.  Furthermore, $P(Y\cup \Gamma)=\{ f\in C(Y\cup \Gamma): f|Y\in P(Y)\}$.
\end{theorem}

\blem\cite[Lemma~2.1]{Izzo:2022}\label{prelemma}
Let $\lambda$ be a closed set in $\R^N$, $N\geq 3$, of topological dimension at most $1$, let $a$ and $a'$ are two points in $\lambda$, and let $\Omega$ be a connected open set of $\R^N$ that contains $a$ and $a'$.  Then there is an arc $J$ from $a$ to $a'$ contained in $\Omega$ that intersects $\lambda$ only in the end points $a$ and $a'$ of $J$ and is such that the open arc $J\sm \{a,a'\}$ is $C^\infty$-smooth.
\elem

The following elementary lemma will also be needed.

\blem\label{elem-approx-lemma}
Let $X\subset\C^N$ be a compact set such that $P(X)=C(X)$, and let $\row xk$ be points in $\C^N$.  Then $P(X\cup\{\row xk\})=C(X\cup\{\row xk\})$.
\elem

\bpf
By induction it suffices to consider the case $k=1$.  Assume that $x_1\notin X$ since otherwise there is nothing prove.  Since $P(X)=C(X)$, the set $X$ is polynomially convex, so there exists a polynomial $p$ such that $p(x_1)=1>\|p\|_X$.  The function $f$ that is 1 at $x_1$ and identically 0 on $X$ is the uniform limit of the sequence $(p^n)_{n=1}^\infty$ and hence is in $P(X\cup\{x_1\})$.  

Given $g\in C(X\cup\{x_1\})$, choose a sequence of polynomials $(g_n)$ such that $g_n\rightarrow g$ uniformly on $X$.  Then $g_n(1-f) + g(x_1)f\rightarrow g$ uniformly on $X\cup\{x_1\}$, so $g$ is in $P(X\cup\{x_1\})$.
\epf

\bpf[Proof of Theorem~\ref{rectifiable-approximated}]
We may assume without loss of generality that the diameter of $B$ is less than $\vep$.
The simple closed curve $\gamma_a$ will be obtained from $\gamma$ by removing an  arc from $\gamma$ contained in $B$ and replacing it with a new arc having the same end points.  It is then immediate that $\gamma_a\sm B=\gamma \sm B$
and that under a suitable parametrization $\|\gamma-\gamma_a\|_\infty<\vep$.  We must show that the removal and insertion can be done in such a way that the other conditions are satisfied.

Remove from $\gamma$ an open arc lying in $B$ to obtain a rectifiable arc $\lambda$ whose end points we will denote by $p$ and $q$.  As noted immediately after Corollary~\ref{arc-union}, there exists an arc $\tau$, which can be taken to lie in $B$ and disjoint from $\lambda$ such that $\tau$ is not rationally convex and such that, denoting the end points of $\tau$ by $a$ and $d$, there are distinct points $b$ and $c$ in the interior of $\tau$ such that traversing $\tau$ from $a$ to $d$, one encounters $b$ before $c$ and the arcs $ac$ and $bd$ are polynomially convex and satisfy $P(ac)=C(ac)$ and $P(bd)=C(bd)$.

The set $B \sm (\lambda \cup \tau)$ is connected (because a connected manifold of real dimension greater than or equal to three cannot be disconnected by a subspace of topological dimension one \cite[Corollary~1,~p.~48]{Hurewicz-Wallman:1948}).  
Therefore, by repeated application of Lemma~\ref{prelemma}, we can obtain smooth arcs $\ell_1$ and $\ell_2$ in $B$ such that the arc $\ell_1$ has end points $p$ and $a$, the arc $\ell_2$ has end points $q$ and $d$, and the set
\begin{equation*}
\gamma_a= \lambda\cup\ell_1\cup\tau\cup\ell_2
\end{equation*}
is a simple closed curve.

Since $\tau$ is not rationally convex, $\gamma_a$ is not rationally convex by Lemma~\ref{convexity-of-subsets}.

Let $\gamma_1=\lambda\cup\ell_1\cup\ell_2\cup ac$ and $\gamma_2=\lambda\cup\ell_1\cup\ell_2\cup bd$.  Then the interiors of $\gamma_1$ and $\gamma_2$ cover $\gamma_a$.  Set $Y=ac\cup \{d,p,q\}$ and $\Gamma=\lambda\cup\ell_1\cup\ell_2$.  Then $Y$ is polynomially convex and $P(Y)=C(Y)$ by Lemma~\ref{elem-approx-lemma}, so Theorem~\ref{from-Stout-paper} applies and yields that $\gamma_1$ is polynomially convex and satisfies 
$P(\gamma_1)=C(\gamma_1)$.  Similarly, $\gamma_2$ is polynomially convex and satisfies $P(\gamma_2)=C(\gamma_2)$.
\epf

%
%

\section{A nonrationally convex ball and chain}\label{ball-n-chain}

In this section we consider polynomial convexity of the union of a polynomially convex arc and the closed unit ball $\ob$.  We begin with an observation of Stout.

\bthm\label{Stout-observation}
Let $J$ be an arc in $\C^N$ that meets $\ob$ in a single end point, say $a$.  If the half-open arc $J\sm\{a\}$ is locally of finite length, then $\ob\cup J$ is polynomially convex.
\ethm

Note that the hypotheses do not require that $J$ itself has finite length.

\bpf
The cohomology group $\check H^1(\ob\cup J; Z)$ is the zero group, so the
natural map $\check H^1(\ob\cup J;Z) \rightarrow \check H^1(\ob;Z)$ is a monomorphism, whence application of Theorem~\ref{from-Stout-paper}
yields the result.
\epf

Variations on Theorem~\ref{Stout-observation} can be proved similarly.  For instance, it 
is enough to require only that the interior of $J$ be locally of finite length.  Also the arc $J$ could be taken to be disjoint from $\ob$ or could intersect $\ob$ in an interior point of $J$ with the two half-open arcs making up $J$ each required to be locally of finite length.

We turn now to the counterexamples.

\bthm\label{ball-disjoint-arc}
There exists a polynomially convex arc $J$ in $\C^3$ that is disjoint from the closed unit ball $\ob$ such that $\ob\cup J$ is not rationally convex.
\ethm

Since there exists an arc in the plane having positive planar measure (a theorem of Osgood \cite{Osgood:1903}), the above result is a special case of the following one.

\bthm
Let $\Gamma$ be a compact subset of $\C$ of positive planar measure with empty interior and connected complement.  There exists a polynomially convex set $E$ homeomorphic to $\Gamma$ in $\C^3$ and disjoint from the closed unit ball such that $\ob\cup E$ is not rationally convex.
\ethm

\bpf
Set
$$g(\zeta)=\int\!\!\!\int_\Gamma\, \frac{dx\,dy}{z-\zeta}.$$
Then $g$ is a nonconstant, continuous function on the Riemann sphere $S^2$ that is holomorphic off $\Gamma$ and vanishes at $\infty$.  (See Section~2.)  Choose a point $z_0\in \C\sm\Gamma$ such that $g(z_0)\neq0$.  Set
$$h(z)=\frac{g(z)-g(z_0)}{z-z_0}.$$
Then $h$ is also a nonconstant, continuous function on the Riemann sphere $S^2$ that is holomorphic off $\Gamma$ and vanishes at $\infty$.  Note that $\infty$ is the only common zero of $g$ and $h$.  In particular, $g$ and $h$ have no common zero on $\Gamma$, so there is a positive constant $C$ such that $C|(g,h)|>1$ everywhere on $\Gamma$.

Since $(g,h)(\infty)=0$, we have for sufficiently large $R>0$ that $|(Cg,Ch)|<1/2$ everywhere on the circle $\{|z|=R\}$.  Fix $R>0$ large enough that, in addition, the circle $\{|z|=R\}$ encloses the set $\Gamma$.

Define $F:\C\to\C^3$ by
$$F(z)=\bigl(z/2R, Cg(z), Ch(z)\bigr),$$
and set $E=F(\Gamma)$.
Then $F$ maps $\Gamma$ homeomorphicly onto $E$.  Furthermore, $E$ is polynomially convex because polynomials in $z_1$ are dense in $C(E)$, since polynomials in $z$ are dense in $C(\Gamma)$ by Lavrentiev's theorem \cite[Theorem~II.8.7]{Gamelin:1984}.
In addition, because $|(Cg,Ch)|>1$ everywhere on $\Gamma$, the set $E$ is disjoint from $\ob$.

Note that $F(\{|z|=R\})$ is contained in $\ob$. 
Thus the connected set $F(\{|z|\leq R\})$ intersects each of the two disjoint closed sets $\ob$ and $E$, and hence $F(\{|z|\leq R\})$ cannot be contained in the union of $\ob$ and $E$.  Thus to show that $\ob\cup E$ is not rationally convex, it suffices to show that $F(\{|z|\leq R\})$ lies in the rational hull of $\ob\cup E$.

Suppose $p$ is a polynomial on $\C^3$ having no zeros on $\ob$ and $E$.  Since each of $\ob$ and $E$ is simply coconnected, $p$ has a continuous logarithm in a neighborhood of $\ob\cup E$.  Consequently, $p\circ F$ has a continuous logarithm on a neighborhood of $\{|z|=R\}\cup \Gamma$ in $\C$.  By the argument principle it follows that $p\circ F$ has no zeros on $\{|z|\leq R\}$.  Thus $p$ has no zeros on $F(\{|z|\leq R\})$. Consequently, $F(\{|z|\leq R\}$) is contained in the rational hull of $\ob\cup E$.
\epf

\bcor\label{non-rat-ball-n-chain}
There exists a polynomially convex arc $J$ in $\C^3$ that meets the closed unit ball $\ob$ in a single end point such that $\ob\cup J$ is not rationally convex.
\ecor

\bpf
By Theorem~\ref{ball-disjoint-arc} there is a polynomially convex arc $E$ in $\C^3$ disjoint from $\ob$ such that $\ob\cup E$ is not rationally convex.  One can construct an arc $\sigma$ from a point of $\pb$ to an end point of $E$ such that the interior of $\sigma$ is disjoint from $\ob\cup E$ and $\sigma$ is smooth except possibly at the end point where it meets $E$.
The details of the construction are similar to the proof of 
\cite[Theorem~1.2]{Izzo-Stout:2022b} so we omit them; the basic idea is to choose a sequence of points in the complement of $\ob\cup E$ converging to an end point of $E$, and then choose smooth arcs connecting successive points of the sequence and fitting together so as to form a smooth arc.

Set $J=E\cup\sigma$.  Then $J$ is an arc and is polynomially convex by \cite[Theorem~1.7]{Izzo-Stout:2022b} (a corollary of Theorem~\ref{from-Stout-paper}).  Of course 
$\hr{\ob\cup J} \supset\hr{\ob\cup E}$, and the later set cannot be contained in $(\ob\cup E) \cup \sigma= \ob \cup J$ by Theorem~\ref{gen-Rossi-consequence}.
Thus $\hr{\ob\cup J}\sm (\ob\cup J)$ is nonempty.
\epf

%
%

\section{Unions with hull without analytic discs}\label{no-discs}

\bthm\label{no-discs3}
There exist two polynomially convex Cantor sets $K_1$ and $K_2$ in $\C^3$ such that the polynomial hull of $K_1\cup K_2$ is nontrivial but contains no analytic discs.
\ethm

In $\C^4$ the two polynomially convex Cantor sets can be taken to be disjoint.

\bthm\label{no-discs4}
There exist two disjoint polynomially convex Cantor sets $K_1$ and $K_2$ in $\C^4$ such that the polynomial hull of $K_1\cup K_2$ is nontrivial but contains no analytic discs.
\ethm

In connection with these results, note that the union of two Cantor sets is itself a Cantor set.

As a corollary of Theorem~\ref{no-discs4} we will prove a general result about unions with polynomial hull without analytic discs.

\bcor\label{no-discs-gen}
Let $(J,K_1,K_2)$ be a triple with $J$ a compact subspace of $\C^N$ and $K_1$ and $K_2$ closed subspaces of $J$ whose union is $J$ and such that each of $J\sm K_1$ and $J\sm K_2$ is uncountable.  Then there exists an embedding $\pi$ of $J$ into $\C^{N+5}$ such that the polynomial hull of $\pi(J)$ is not trivial but contains no analytic discs while each of $\pi(K_1)$ and $\pi(K_2)$ is polynomially convex and $P(\pi(K_j))=C(\pi(K_j))$, $j=1,2$. 
\ecor

The proofs of these results are based on a result from the paper of the author and Norman Levenberg \cite{Izzo-Levenberg:2019} which associates to each compact set with nontrivial polynomial hull in $\C^N$ a compact set in $\C^{N+1}$ with a nontrivial polynomial hull that contains no analytic discs.  We quote the needed result here for the reader's convenience.

\begin{theorem}[\cite{Izzo-Levenberg:2019}, Theorem~1.1]\label{Wermergen} 
Let $X \subset \C^N$ be a compact set whose polynomial hull is nontrivial.  Then there exists a compact set $Y \subset \C^{N+1}$ such that, letting $\pi$ denote the restriction to $\h Y$ of the projection $\C^{N+1} \to \C^N$ onto the first $N$ coordinates, the following conditions hold:
\item {\rm (i)} $\pi (Y) = X$
\item {\rm (ii)} $\pi (\hh Y) = \hh X$
\item {\rm (iii)} $\h Y$ contains no analytic discs
\item {\rm (iv)} each fiber $\pi^{-1} (z)$ for $z \in \h X$ is totally disconnected.
\end {theorem}

We will also need two lemmas concerning perfect subsets.  Recall that a subset of a space is said to be \emph{perfect} if it is closed and has no isolated points.  Every space contains a unique largest perfect subset (which can be empty), namely the closure of the union of all perfect subsets of the space.

\blem\label{trivial-from-Oka-Weil}
Let $Y$ be a polynomially convex set in $\C^N$, and let $K$ be the largest perfect subset of $Y$.  Then $K$ is polynomially convex.
\elem

\bpf
Note that $\h K\subset \h Y=Y$.  Thus since $K$ is the \emph{largest} perfect subset of $Y$, if $\h K$ were strictly larger than $Y$, then $\h K$ would have an isolated point, and the isolated point would necessarily be a point of $\h K\sm K$.  But this is impossible, for it follows from the Oka-Weil theorem that every component of $\h K$ must intersect $K$.
\epf

\blem\label{prove-using-Zorn}
Let $Y$ be a space, and let $Y_1$ and $Y_2$ be closed subspaces of $Y$ whose union is $Y$.  Let $K$, $K_1$, and $K_2$ be the largest perfect subsets of $Y$, $Y_1$, and $Y_2$, respectively.  Then $K_1\cup K_2=K$.
\elem

In case $Y_1$ and $Y_2$ are disjoint, this lemma is essentially obvious.  Only that case is needed for the proof of Theorem~\ref{no-discs4}, but the general case is needed for the proof of Theorem~\ref{no-discs3}.

\bpf
The proof is an application of Zorn's Lemma.  Let $\P$ be the collection of ordered pairs $\vector\Sigma 2$ of closed subsets of $Y$ such that $K_j\subset \Sigma_j\subset Y_j$ for $j=1,2$ and $\Sigma_1\cup\Sigma_2=K$.  Note that $\P$ is nonempty since $(Y_1\cap K, Y_2\cap K)$ is in $\P$.  Partially order $\P$ by declaring $\vector \Sigma 2 \leq \vector {\Sigma'}2$ if $\Sigma_1\subset\Sigma'$ and $\Sigma_2\subset \Sigma'_2$.  Given a totally ordered subcollection $\Q$ of $\P$, set $$Q_1=\bigcap\nolimits_{\vector\Gamma2\in \Q}\Gamma_1\quad {\rm and}\quad Q_2=\bigcap\nolimits_{\vector\Gamma2\in \Q}\Gamma_2.$$
It is easily verified that $\vector Q2$ is a lower bound for $\Q$ in $\P$.  Thus by Zorn's lemma, $\P$ has a minimal element $\vector P2$.  To conclude the proof, it suffices to show that $\vector P2=\vector K2$ since then $K_1\cup K_2=P_1\cup P_2=K$.

Assume to get a contradiction that $\vector P2\neq \vector K2$.  Then without loss of generality, $P_1$ properly contains $K_1$.  Then $P_1$ is not perfect and hence has an isolated point $x$.  Since $K$ is perfect, and $P_1\cup P_2=K$, it must be that $x$ is a limit point of $P_2$, and hence $x$ is in $P_2$.  Consequently, the closed sets $P_1\sm\{x\}$ and $P_2$ have union $K$.  Also since $x$ is isolated in $P_1$, the point $x$ can not be in $K_1$.  Thus $K_1\subset P_1\sm\{x\}$.  We conclude that $(P_1\sm\{x\}, P_2)$ is a member of $\P$ that is strictly smaller than $\vector P2$, contrary to the minimality of $\vector P2$.
\epf

\bpf[Proof of Theorems~\ref{no-discs3} and~\ref{no-discs4}]
Set $N=2$ for the proof of Theorem~\ref{no-discs3}, or $N=3$ for the proof of Theorem~\ref{no-discs4}.  By Corollary~\ref{disjoint-Cantor-sets} and Theorem~\ref{two-Cantor-sets} there exist polynomially convex Cantor sets $X_1$ and $X_2$ in $\C^N$ such that their union, which we will denote by $X$, has nontrivial polynomial hull.  Furthermore, when $N=3$, we can choose $X_1$ and $X_2$ to be disjoint.  Let $Y$ be the set in $\C^{N+1}$ obtained from $X$ by applying Theorem~\ref{Wermergen}.  As in Theorem~\ref{Wermergen} let $\pi:\h Y \to \C^N$ be the restriction to $\h Y$ of the projection $\C^{N+1}\to \C^N$ onto the first $N$ coordinates.  Set $Y_1=\pi^{-1}(X_1)$ and $Y_2=\pi^{-1}(X_2)$.  Let $K$, $K_1$, and $K_2$ be the largest perfect subsets of $Y$, $Y_1$, and $Y_2$, respectively. Then $K=K_1\cup K_2$ by Lemma~\ref{prove-using-Zorn}.  Also $\hh K\supset \hh Y$ by \cite[Lemma~4.2]{Izzo:2020}, so condition (ii) of Theorem~\ref{Wermergen} gives that $\h K$ is nontrivial.  Since $\h K\subset \h Y$, condition (iii) gives that $\h K$ contains no analytic discs.  It is easily verified that the polynomial convexity of $X_1$ and $X_2$ implies that $Y_1$ and $Y_2$ are polynomially convex.  Consequently, $K_1$ and $K_2$ are polynomially convex by Lemma~\ref{trivial-from-Oka-Weil}.  Finally, conditions (i) and (iv) of Theorem~\ref{Wermergen} imply that $Y$ is totally disconnected, and hence the same is true of the subsets $K_1$ and $K_2$.  Consequently, $K_1$ and $K_2$ are Cantor sets by the usual characterization.
Note that $K_1$ and $K_2$ are disjoint when $N=3$.
\epf

\bpf[Proof of Corollary~\ref{no-discs-gen}]
As in the proof of Corollary~\ref{intrinsic-form}, we can choose Cantor sets $G_1$ and $G_2$ in $J\sm K_1$ and $J\sm K_2$, respectively.  Let $\row f4$ be the components of a homeomorphism of $G_1\cup G_2$ onto the nonpolynomially convex set in Theorem~\ref{no-discs4} that maps $G_1$ onto the first Cantor set in that theorem and $G_2$ onto the second Cantors set there.  

Extend $\row f4$ to continuous complex-valued functions $\row{\tilde f}4$ on $J$.  Let $x_1,\ldots, x_N$ denote the real coordinate functions of $\R^N$.  Choose a continuous real-valued function $\rho$ on $J$ whose zero set is precisely $G_1\cup G_2$.  Let $\pi: J\rightarrow \C^{N+5}$ be the mapping whose components are the functions $\row{\tilde f}4, \rho, \rho x_1,\ldots \rho x_N$.  

Because $\pi(J)$ is contained in $\C^4\times \R^{N+1}$, the polynomial hull $\h {\pi(J)}$ is the union of the polynomial hulls of the slices $\pi(J)\cap (\C^4\times \{r\})$ for $r\in \R^{N+1}$ (by \cite[Proposition~3.1]{ISW:2016} for instance).  Each of these slices is a single point with the exception of the slice given by $r=0$.  That slice is the set $\pi(G_1\cup G_2)\cap (\C^4\times \{0\})$, which is exactly the image, under the canonical embedding of $\C^4$ into $\C^{N+5}$, of the set in Theorem~\ref{no-discs4} with nontrivial hull containing no analytic discs.  Consequently, $\pi(J)$ has nontrivial hull containing no analytic discs.

Application of the Bishop antisymmetric decomposition (\cite[Theorem~2.7.5]{Browder:1969} or \cite[Theorem~13.1]{Gamelin:1984}) shows that $P(\pi(K_j))=C(\pi(K_j))$, $j=1,2$, and hence each of $\pi(K_1)$ and $\pi(K_2)$ is polynomially convex.
\epf

\medskip
\centerline{\textsc{Acknowledgment}}
\smallskip
This research was begun while the author was a visitor at the University of Michigan.  He thanks the Department of Mathematics for its hospitality.
He also thanks Lee Stout for inspiring correspondence related to the paper.

\bibliography{Izzo.Nonrationally}

\end{document}